\documentstyle[12pt,twoside]{article}
\newtheorem{theorem}{\bf Theorem}[section]
\newtheorem{proposition}[theorem]{\bf Proposition}

\newtheorem{corollary}[theorem]{\bf Corollary}

\newtheorem{remark}[theorem]{\bf Remark}

\input{amssym}

\date{}
\begin{document}
\title{{\Large\bf Jordan derivations on certain Banach algebras}}

\author{{\normalsize\sc M. J. Mehdipour\footnote{Corresponding author}~, GH. R. Moghimi and N. Salkhordeh}}
\maketitle

{\footnotesize {\bf Abstract.} In this paper, we study the types of Jordan derivations of a Banach algebra $A$ with a right identity $e$. We show that if $eA$ is commutative and semisimple, then every Jordan derivation of $ A $ is a derivation. In this case, Jordan derivations map $A$ into the radical of $A$.  We also prove that every Jordan triple left (right) derivation of $ A $  is a Jordan left (right) derivation. Finally, we investigate the range of  Jordan left derivations and establish that every Jordan left derivation of $ A $ maps $ A $ into $eA$.

}
{}

\section{\normalsize\bf Introduction}

Let $A$ be a Banach algebra with the right annihilator $\hbox{ran}(A)$ and the radical $\hbox{rad}(A)$. Let us recall that a linear map $ d $ on $ A $  is called a \emph{derivation} if for every $a, x\in A$
$$
d(ax)=d(a)x+ad(x);
$$
also, $d$ is a $\textit{Jordan derivation}$ if
$$ d(a^{2})= d(a)a+ad(a),$$
and $ d $ is called a  $\textit{Jordan left derivation}$ if
$$ d(a^{2}) = 2ad(a)$$
for all $a\in A $. Finally, $d$ is called a \emph{Jordan triple left derivation} if
$$
d(a^{3}) =3a^2d(a)
$$
for all $a\in A$.
 Jordan right derivations and Jordan triple right derivations are defined similarly.

Herstein [7] proved that any Jordan derivation of prime rings with a suitable characteristic is a derivation. Cusack [5] extended this result to semiprime rings. Sinclair [13] proved that every continuous Jordan derivation of a semisimple Banach algebra is a derivation; for derivations and Jordan derivations on group algebras see [1, 2, 12].
Jordan left derivations have been introduced and studied by Bre\v{s}ar and Vukman [4]. They proved that there is no nonzero Jordan left derivation of noncommutative prime rings with a suitable characteristic.

In this paper, we consider a Banach algebra $ A $ with a right identity $e$ and study the types of Jordan derivations of it. In Section 2, we show that $ \hbox{ran}(A) $ is invariant under  Jordan derivations of $A$, and prove that if $eA $ is commutative and semisimple, then every Jordan derivation of $ A $ is a derivation with the range into $ \hbox{ran}(A) $. In Section 3, we prove that every Jordan triple left (right) derivation of $ A $ is a Jordan left (right) derivation. We also verify that every Jordan left derivation of $ A $ maps $ A $ into $ eA $.

\section{\normalsize\bf Jordan derivations on Banach algebras}

We commence this section with the following result.

\begin{proposition}\label{s1}
\ Let $A$ be a Banach algebra with a right identity, and let $ d: A \rightarrow A $ be a Jordan derivation.  Then the following statements hold.

\emph{(i)} $ \emph{ran}(A) $ is invariant under $ d $.

\emph{(ii)} If $ d$ maps $ A $ into  $ \emph{ran}(A)  $, then $ d $ is a derivation.
\end{proposition}
{\it Proof.}
(i) Let $ d $ be a Jordan derivation of $ A $. Then for every $ a, x\in A $, we have
\begin{equation}\label{m1}
 d(ax+ xa) = d(a) x + a d(x) + d(x) a+ x d(a).
\end{equation}
Let $e$ be a right identity of $ A $. Putting $ a=x=e $ in $(1)$, we have
\begin{eqnarray}
2 d(e) &=& d(ee + ee ) \nonumber \\
&=& d(e)e + ed(e)+ d(e)e + ed(e)   \nonumber \\
&=& 2 d(e) + 2 ed(e).  \nonumber
\end{eqnarray}
Thus $ed(e)= 0$. Hence $  d(e) \in  \hbox{ran}(A) $. Taking $ a=e $ in $(1)$, we obtain
\begin{eqnarray}
d(ex)+ d(x)&=& d(ex+ xu) \nonumber \\
&=& d(e) x + ed(x) + d(x)e + x d(e)  \nonumber \\
&=& d(e) x + ed(x) + d(x)  \nonumber
\end{eqnarray}
for all $ x \in A $. So
\begin{equation}\label{m2}
d(ex)= d(e) x + ed(x).
\end{equation}
Consequently, for every $ r \in  \hbox{ran}(A) $ and $ a\in A $, we have
\begin{eqnarray}
ad(r) &=& aed(r) \nonumber \\
&=& a(d(e) r + ed(r))  \nonumber \\
&=& ad(er) \nonumber \\
&=& 0. \nonumber
\end{eqnarray}
This shows that $ d(r) \in  \hbox{ran}(A) $.\\
(ii) Let $ d $ map $ A $ into $ \hbox{ran}(A) $. Then by $(1)$ and $(2)$ for every $ a\in A $ and $ r \in  \hbox{ran}(A)$
$$ d(ea) = d(e)a\quad \hbox{and}\quad d(ra) = d(r)a.  $$
Note that for every $ a\in A $, there exists $ r \in  \hbox{ran}(A) $ such that $ a = ea + r $. Hence for every $ x \in A $, we have
\begin{eqnarray}
d(ax)&=& d(eax+r x ) \nonumber \\
&=& d(eax) + d (r x )  \nonumber \\
&=& d(e)a x + d(r) x  \nonumber \\
&=& d(ea) x + d(r) x \nonumber \\
&=& d(ea + r) x \nonumber \\
&=& d(a) x. \nonumber
\end{eqnarray}
Since $ ad (x) = 0 $, it follows that $ d(ax) = d(a) x + ad(x) $. Hence, $ d $ is a derivation of $ A $.
$\hfill\square$\\

A linear mapping $T$ on a Banach algebra $A$ is called $\textit{ spectrally bounded }$ if there exists a non-negative number $\alpha$ such that
$$r(T(a)) \leq \alpha r(a),$$
for all $a\in A$, where $r(.)$ denotes the spectral radius.

\begin{theorem}\label{s2}
Let $ A $ be Banach algebra with a right identity $e$. If $eA $ is commutative and semisimple, then every Jordan derivation $ d $ on $ A $ is a derivation. Furthermore, $ d $ is spectrally infinitesimal and $ d(A) $ is contained in $\emph{ran}(A) $.
\end{theorem}
{\it Proof.}
Let $ d $ be a Jordan derivation of $ A $. In view of Theorem 2.1, $d $ maps $\hbox{ran}(A) $ into $\hbox{ran}(A) $. Hence the Jordan derivation $  D: A/ \hbox{ran}(A) \rightarrow A/ \hbox{ran}(A)  $ defined by $$ D(a+ \hbox{ran}(A)) = d(a) + \hbox{ran}(A)$$ is well-defined. Note that $ A/ \hbox{ran}(A)$ and $eA$ are isomorphic as Banach algebras. This together with Corollary 1.5.3 (ii) of [6] implies that $ A/ \hbox{ran}(A) $ is a commutative semisimple Banach algebra. Hence $ D $ is a derivation and so $ D $ is zero on $ A/ \hbox{ran}(A) $; see [3, 14, 15]. This implies that $ d(a) \in \hbox{ran}(A) $ and by Proposition 2.1, $ d $ is a derivation. To complete the proof, we only note that $ d(A) $ is nilpotent and so $ d $ is spectrally infinitesimal.
$\hfill\square$\\

For a Banach algebra $A$ and a positive integer $k$, a linear mapping $T$ on $A$ is called \textit{ $k-$centralizing }  if
$$T(a)a^{k}-a^kT(a)\in Z(A)$$ for all $a\in A$, where $Z(A)$ is the center of $A$.

\begin{corollary}\label{s3}
Let $ A $ be Banach algebra with a right identity $e$. If $ eA $ is commutative and semisimple, then the following statements hold.

\emph{(i)} The range of a Jordan derivation of $ A $ is contained in $ \emph{rad}(A) $.

\emph{(ii)} If $d_{1}$ and $d_{2}$ are two Jordan derivations on $A$, then $d_{1}d_{2}$ is a derivation.

\emph{(iii)} The zero map is the only $k-$centralizing Jordan derivation of $A$, where $ k $ is a positive integer.
\end{corollary}
{\it Proof.} The statement (i) follows at once from Theorem 2.2. For (ii), note that if $d_1$ and $d_2$ are two Jordan derivations on $A$, then they are derivations on $A$ with the range into $\hbox{ran}(A)$ by Theorem 2.2. Thus for every $a, b\in A$, we have
$$
d_{1}(a)d_{2}(x)+d_{1}(x)d_{2}(a) = 0.
$$
This implies that $d_{1}d_{2} $ is a derivation of $A$. Finally, assume that $d$ is a $k-$centralizing Jordan derivation of $A$, where $ k $ is a positive integer. Then by Theorem 2.2, for every $a\in A$ we have
\begin{eqnarray*}
d(a)a^k&=&d(a)a^k-a^kd(a)\\
&\in&\hbox{ran}(A)\cap Z(A)\\
&=&\{0\}.
\end{eqnarray*}
Thus $d(e)=0$. Let $a\in A$ and set $r=a-ea$. Then  $$r+e=(r+e)^k.$$ This implies that $d(r)=0$. So
$$
d(a)=d(r+ea)=d(ea)=d(e)a=0.
$$
Therefore, $d=0$. That is, (iii) holds.
$\hfill\square$

\begin{theorem}\label{s4}
Let $ A $ be Banach algebra such that $ A/ \emph{ran}(A) $ is commutative and semisimple. Then the following statements hold.

\emph{(i)} Every derivation of $ A $ maps $ A $ into $ \emph{ran}(A) \subseteq \emph{rad}(A) $.

\emph{(ii)} Every derivation of $ A $ is spectrally infinitesimal.

\emph{(iii)} The product of two derivations of $A$ is a derivation of $A$.

\end{theorem}
{\it Proof.}
A simple computation shows that $ \hbox{ran}(A)$ is invariant under derivations on $ A $. Hence we can apply the proof of Theorem 2.2 to derivations instead of Jordan derivations. Therefore, (i) and (ii) hold.
The statement (iii) follows from (i).
$\hfill\square$\\

Let $A$ be a Banach algebra and let $k$ be a positive integer. A linear mapping $T$ on $A$ is called $\textit{ $k-$skew centralizing }$  if
$$T(a)a^{k}+a^kT(a) \in Z(A)$$ for all $a\in A$, in a special case, if for every $a\in A$
$$T(a)a^{k}+a^kT(a)=0,$$
then T is called $k$-\emph{skew commuting}. Also, $T$ is called \emph{spectrally infinitesimal} if $r(T(a)) = 0$ for all $a\in A$.

\begin{remark}{\rm  In Theorems 2.2, 2.4 and Corollary 2.3, we consider a Banach algebra $ A $ such that either\\

 $eA $ is semisimple when $ A $ has a right identity $e$ or $ A/ \hbox{ran}(A) $ is semisimple. ($\maltese$)\\\\
 Let us remark that if $ G $ is a locally compact group, then the Banach algebras $ L^{1}(G) $, $M(G)$ and $ L_{0}^{\infty}(G)^{*}$ satisfy the condition ($\maltese$); For an extensive study of these Banach algebras see [8, 9, 10, 11].
Also, $ C^{*} $-algebras, semisimple Banach algebras with a bounded approximate identity and Banach algebras with  $$ \hbox{rad}(A) =  \hbox{ran}(A) $$ satisfy the condition ($\maltese$). Under this condition,
 one can prove the results of [1] for $ A $ instead of Banach algebra $ L_{0}^{\infty} (G)^{*}$. In the following, we state some of the important results.

(a) Every generalized derivation of $ A $ is spectrally bounded.

(b) A generalized derivation $ (D,d) $ of $ A $ is spectrally infinitesimal if and only if $ D=d $; or equivalently, $ D(A) \subseteq \hbox{rad}(A)$.

(c) A generalized derivation $ (D,d) $ of $ A $ is $k-$centerlizing if and only if $ D = R_{x} $ for some $ x \in Z(A) $, where $R_x(a)=ax$ for all $a\in A$ and $ k $ is a positive integer.

(d) If a generalized derivation $ (D,d) $ of $ A $ is $k-$skew centralizing, then $ D = R_{x} $ for some $ x \in Z(A) $.

(e) If a generalized derivation $ (D,d) $ of $ A $ is $k-$skew commuting, then $ D = 0 $.}
\end{remark}

\section{\normalsize\bf Jordan triple left (right) derivations of Banach algebras}

We first investigate Jordan left derivations.

\begin{theorem}\label{n1}
Let $A$ be a Banach algebra with a right identity $e$, and let $d$ be a Jordan left derivation of $ A $. Then $ d $ maps $ A $ into $ eA $.
\end{theorem}
{\it Proof.}
 Let $d$ be a Jordan left derivation of $ A $. Then
\begin{equation}\label{m13}
d(e) =2ed(e)
\end{equation}
and so
$$ ad(e) =2ad(e) $$
for all $a\in A$. It follows that
$$d(e) \in \hbox{ran}(A) $$
and hence $d(e) =0 $ according to $(3)$. Now, fix $r \in \hbox{ran}(A)$. Then
\begin{eqnarray}
d(r)&=& d(r)+d(e) \nonumber \\
&=& d((r +e)^{2})  \nonumber \\
&=& 2(r +e)·d(r +e) \nonumber \\
&=& 2r ·d(r)+2ed(r). \nonumber
\end{eqnarray}
So for every $a\in A$,
$$ ad(r) = 2ad(r),$$
which shows that
$$ d(r) \in \hbox{ran}(A) $$
for all $r \in \hbox{ran}(A) $. By hypothesis, we have
\begin{equation}\label{m14}
d(a^{2}) =2ad(a),
\end{equation}
for all $a\in A$. Let us replace $a$ by $a+e$ in $(4) $. Then
\begin{eqnarray}
d(a+ea) &=& 2ad(e)+2ed(a) \nonumber \\
&=& 2ed(a). \nonumber
\end{eqnarray}
 Hence
\begin{equation}\label{m15}
d(a) =2ed(a)-d(ea),
\end{equation}
for all $ a\in A$. From the fact that
$$a-ea\in \hbox{ran}(A),$$
we infer that
\begin{equation}\label{m16}
ed(a-ea) =0.
\end{equation}
Using $(5) $ and $(6) $, we arrive at
$$d(a-ea) = 2ed(a- ea)- d(e·(a-ea)) = 0.$$
Hence
$$d(a) =d(ea),$$
which together with $(5) $ shows that
\begin{equation}\label{m17}
d(a) =e·d(a),
\end{equation}
for all $a\in A$.
$\hfill\square$

\begin{corollary}\label{n2}
Let $A$ be a Banach algebra with a right identity $e$, and let $d$ be a Jordan left derivation of $ A $. If $ d $ maps $ A $ into $ \emph{ran}(A)$, then $ d $ is zero.
\end{corollary}
{\it Proof.}
Let $ d(A) \subseteq\hbox{ran}(A)$.  Then by Theorem 3.1, we have
$$
d(A) \subseteq\hbox{ran}(A)\cap eA= \lbrace 0 \rbrace,
$$
as claimed.
$\hfill\square$

\begin{theorem}\label{n4}
Let $A$ be a Banach algebra with a right identity. If $d$ is a Jordan triple left derivation of $A$, then $ d $ is a Jordan left derivation.
\end{theorem}
{\it Proof.} Let $d$ be a Jordan triple left derivation of $A$. Then
for every $a\in A$, we have
\begin{equation}\label{m18}
d(a^{3}) =3a^{2} d(a).
\end{equation}
Replace $a$ by $a+e$ in $(8)$. Then
\begin{equation}\label{m19}
2d(a^{2})+d(a)+2d(ea)+d(e·a^{2}) =3ed(a)+3ad(a)+3e·ad(a).
\end{equation}
Note that $d(e) =0$. The substitution $-a$ for $ a$ in $(8)$ gives
\begin{equation}\label{m20}
2d(a^{2})-d(a)-2d(e·a)+d(e·a^{2}) = 3ed(a)+3ad(a)+3e·ad(a).
\end{equation}
Combining $(9)$ with$(10)$ we arrive at
\begin{equation}\label{m21}
d(a)+2d(ea) =3ed(a).
\end{equation}
If we take $a-ea$ instead of $a$ in $(11)$, then
\begin{equation}\label{m22}
d(a) =d(ea).
\end{equation}
From this and $(11)$ we infer that
\begin{equation}\label{m23}
d(a) =ed(a).
\end{equation}
Now, put $ a=e $ in $ (8) $. Then $ed(e)= d(e)=0$. Thus
\begin{eqnarray}
d(a^{3}+2a^{2}+a+ea^{2}+2ea+e^{3}) &=& 3a^{2}d(a)+3a^{2}d(e)+3ed(a)+3e^{2}d(e) \nonumber \\
&+& 3aed(a)+3ad(e)+3ead(a)+3ead(e). \nonumber
\end{eqnarray}
Since $ d(a^{3}) = 3a^{2}d(a) $, it follows that
$$
2d(a^{2}) + d(a) + d(ea^{2})+2d(ea)= 3ed(a)+3ad(a)+3ead(a).
$$
This together with $ (12) $ and $ (13) $ implies that
$$ 3d(a^{2})+3d(a) =  3d(a)+3ad(a)+3ead(a).$$
Hence
\begin{equation}
d(a^{2}) = ad(a)+ ead(a).
\end{equation}
This together with $ (13) $ shows that
 \begin{equation}\label{m25}
d(a^{2})= ed(a^{2})= 2ead(a).
\end{equation}
It follows that
$$
2a^{2}d(a) = a( 2ead(a))= ad(a^{2}).
$$
Thus
$$ 2a^{2}d(a) = ad(a^{2}). $$
Consequently,
$$ 2d(a^{3}) =6a^2d(a)=3(2a^2d(a))= 3 ad(a^{2}).  $$
Hence
\begin{eqnarray}
2d((a+e)^{3}) &=&  3(a+e) d((a+e)^{2})\nonumber \\
&=& 3ad((a+e)^{2})+3ed((a+e)^{2}) \nonumber \\
&=&  3ad((a+e)^{2})+3 d((a+e)^{2}).   \nonumber
\end{eqnarray}
So
\begin{eqnarray}
2 d(a^{3}+2a^{2}+2ea+ea^{2}+a) &=& 3 ad(a^{2}+2a)+3 d(a^{2}+2a)   \nonumber
\end{eqnarray}
This implies that
 $$ 3d(a^{2})= 6ad(a). $$
Therefore,
$ d(a^{2}) = 2ad(a).$
$\hfill\square$\\

We now prove some results about Jordan right derivations.

\begin{theorem}\label{mj1}
 Let $A$ be a Banach algebra with a right identity, and let $ d $ be a linear map on $ A $. If $d$ is a Jordan triple right derivation, then $d$ is a Jordan right derivation.
 \end{theorem}
{\it Proof.}
Let $d$ be a Jordan triple right derivation of $A$. Then $d (e) =0$ and for every $a\in A$
$$ d((a+e)^{3}) =3d(a+e)(a+e)^{2}.$$
This implies that
\begin{equation}\label{m10}
2d(a^{2})+2d(ea)+d(e·a^{2}) =6d(a)a+2d(a).
\end{equation}
Replacing $a$ by $-a$ in $(16)$ yields
\begin{equation}\label{m11}
2d(a^{2})-2d(ea)+d(e·a^{2}) =6d(a)·a-2d(a).
\end{equation}
From $(16)$ and $(17)$ we obtain
$$d(a) =d(ea)$$ and
$$ 2d(a^{2})+d(e·a^{2}) =6d(a)·a.$$
Therefore,
$$ d(a^{2}) =2d(a)·a.$$
That is, $d$ is a Jordan right derivation.
$\hfill\square$

\begin{theorem}\label{mj2}
 Let $A$ be a Banach algebra with a right identity, and let $ d $ be a Jordan right derivation of $A$. Then $ d $ is zero on $ \emph{ran}(A) $.
\end{theorem}
{\it Proof.}
Let $d$ be a right derivation. It is easy to check that $$d(e) =0 .$$ Also, note that if $r \in \hbox{ran}(A)$, then
$$(r +e)^{2} =r +e$$
and hence $$d(r) =d((r +e)^{2}) =2d(r +e)·(r +e) =2d(r).$$
Therefore, $d(r) =0.$
$\hfill\square$

\vspace{2mm}

{\footnotesize
\noindent {\bf Mohammad Javad Mehdipour}\\
Department of Mathematics,\\ Shiraz University of Technology,\\
Shiraz
71555-313, Iran\\ e-mail: mehdipour@sutech.ac.ir\\
\noindent {\bf Gholam Reza Moghimi}\\
Department of Mathematics,\\ Payame Noor University,\\ Tehran 9395-4697, Iran\\ moghimimath@pnu.ac.ir

\noindent {\bf Narjes Salkhordeh}\\
Department of Mathematics,\\ Shiraz University of Technology,\\ Shiraz
71555-313, Iran,\\ n.salkhordeh@sutech.ac.ir

\end{document}